\documentclass[10pt]{article}

\usepackage{graphics}
\usepackage{graphicx}

\usepackage[a4paper, left=35mm,right=35mm,top=34mm,bottom=34mm]{geometry}
\usepackage[utf8]{inputenc}
\usepackage[T1]{fontenc}
\usepackage[english]{babel}

\usepackage{enumerate}
\usepackage{graphicx}
\usepackage{hyperref}
\hypersetup{
    colorlinks=true,
    linkcolor=blue,
    filecolor=magenta,      
    urlcolor=cyan,
}
\usepackage{listings}
\usepackage{color}

\definecolor{dkgreen}{rgb}{0,0.6,0}
\definecolor{gray}{rgb}{0.5,0.5,0.5}
\definecolor{mauve}{rgb}{0.58,0,0.82}
\lstdefinelanguage{MRGC++}{%
  language=C++,
  morekeywords={T, U, MPI_Irecv, MPI_Isend, MPI_Allreduce, MPI_Waitall, Compute, Map, abs, max, Swap, MPI_Recv_init, MPI_Send_init, MPI_Startall, Copy, Init, InitRecv, InitSend, InitAllReduce, Send, Recv, AllReduce, Finalize, InitSnapshot, Snapshot, SwitchAsync, SnapReduce, MPI_Test, MPI_Start}
}
\usepackage{mathtools,amsthm,amssymb,amsfonts}
\usepackage{algorithm}
\usepackage{algorithmic}
\makeatother
\theoremstyle{plain}

\theoremstyle{definition}

\theoremstyle{remark}

\def\c{\ensuremath{\mathbf{c}}}
\def\f{\ensuremath{\mathbf{f}}}
\def\b{\ensuremath{\mathbf{b}}}

\def\r{\ensuremath{\mathbf{r}}}
\def\w{\ensuremath{\mathbf{w}}}
\def\x{\ensuremath{\mathbf{x}}}
\def\y{\ensuremath{\mathbf{y}}}
\def\xi{\ensuremath{\mathbf{x}_{i}}}
\def\xj{\ensuremath{\mathbf{x}_{j}}}
\def\transp{\ensuremath{^\mathrm{T}}}
\def\lamb{\mbox{\boldmath{\ensuremath{\lambda}}}}
\def\A{\ensuremath{\mathbf{A}}}

\def\P{\ensuremath{\mathbf{P}}}
\def\Q{\ensuremath{\mathbf{Q}}}

\usepackage{caption}
\usepackage{fancyhdr}

\author{
  {\normalsize Guillaume Gbikpi-Benissan}\thanks{Ecole Centrale Paris, France
    (correspondence, frederic.magoules@hotmail.com).}
  \and
  {\normalsize Fr\'ed\'eric Magoul\`es}\footnotemark[1]
}
\title{Coarse Space Correction for Graphic Analysis}
\date{}

\begin{document}
\maketitle
\thispagestyle{fancy}

\begin{abstract}
\noindent In this paper we present an effective coarse space correction addressed to accelerate the solution of an algebraic linear system.
The system arises from the formulation of the problem of interpolating scattered data by means of Radial Basis Functions.
Radial Basis Functions are commonly used for interpolating scattered data during the image reconstruction process in graphic analysis.
This requires to solve a linear system of equations for each color component and this process represents the most time-consuming operation. 
Several basis functions like trigonometric, exponential, Gaussian, polynomial are here investigated to construct a suitable coarse space correction to speed-up the solution of the linear system.
Numerical experiments outline the superiority of some functions for the fast iterative solution of the image reconstruction problem.
\end{abstract}


\section{Introduction}

Interpolation of scattered data is a main issue in image reconstruction theory. The use of Radial Basis Functions (RBFs) for this purpose was introduced in~\cite{Savchenko95} and~\cite{Turk-O'Brien99}. From these papers it appears that solving the System of Linear Algebraic Equations (SLAE) induced by this method comes out to be the most time consuming operation of the whole reconstruction process. Indeed, interpolation of an image by RBF involves performing $O(N^3)$ arithmetic operations, where $N$ denotes the number of data points. Therefore, the computation becomes impractical over several thousands of points.

In spite of this extreme computational cost, RBFs have been widely adopted because of the good results they generally provide, even in many other areas~\cite{Duchon77,Hardy90}. Several advances have been made, allowing to address quite larger data sets, like the use of Compactly-Supported Radial Basis Functions (CSRBFs) proposed by Wendland in~\cite{Wendland95}, from which the resulting SLAE becomes sparse. With this new property, Morse et al. carried out the reconstruction of implicit surfaces from sets of several thousands of points~\cite{Morse01}.

While direct methods used in these approaches allowed to afford less than forty thousands points, iterative methods have been successfully applied for even larger data sets~\cite{Carr01,Beatson99,Ohtake03}.
Recently, the main attention has been centered on partition of unity method where small solutions are stickled together as proposed by Wendland~\cite{Wendland02}. Ohtake et al.~\cite{Ohtake03a} have developped multilevel partition of unity implicit. If the density of the points is not uniform, iterative inclusion of new centers is used to estimate the small solutions. Adaptive iterative inclusion has also been proposed by Hon et al.~\cite{Hon03} for solving large RBF collocation problems but its convergence behavior needs improvement by preconditioning techniques. Hybrid iterative-direct methods, such as domain decomposition methods~\cite{SBG1996,QV1999,TW2005,MR2006} have been widely used to solve large scale linear systems. Additional preconditioning techniques based on transmission conditions~\cite{MM2006}--optimized with a continuous approach~\cite{CN1998,MIT2004,Magoulestopping,magoules:journal-auth:28} or with an algebraic approach~\cite{RMSB2005,MRS2006b,magoules:journal-auth:17,magoules:journal-auth:20,magoules:journal-auth:30}--or on coarse space techniques~\cite{Widlund2008,2009arXiv0911.5725M} have shown strong efficiency and robustness. Magoul\`es et al. in \cite{magoules:journal-auth:11,magoules:journal-auth:25} propose an efficient algorithm to solve the SLAE resulting from the formulation of the problem of image reconstruction from scattered data by means of CSRBF; but the authors did not present a suitable choice of coarse space basis.
In this paper we investigate several original coarse space basis functions and compare their respective efficiencies.

The paper is organized as follows. In section~II the formulation of the CSRBF-based interpolation problem is introduced. The coarse space correction is described in section~III together with the iterative method considered in this paper. Various coarse space basis functions are proposed and compared in section~IV. Finally, section~V contains the conclusions.

\section{Compactly Supported Radial Basis Functions}

In a generalized form, the interpolation problem consists in reconstructing a function from a finite set of linear measurements~\cite{Kybic02a,Kybic02b}. This reconstructed function can be obtained by a linear combination of basis functions, such as in~\cite{Strohmer97,Lee97,Vazquez02,Unser02,Ichige03}. The present study considers Compactly-Supported Radial Basis Functions (CSRBFs)~\cite{Wendland95}, represented by the formula
$$
s(\x)= p(\x) + \sum_{i=1}^{N} \lambda{_i} \phi(\|\x-\xi\|),
$$
where $s$ denotes the CSRBF, $p$, a polynomial of degree one, $\phi$, a radially symmetric function (called basis function), $\lambda{_i}$'s, the CSRBF coefficients, $\xi$'s, the centers of the basis function and the symbol $\|\cdot\|$, the Euclidean norm of a vector.
Defining an interpolating CSRBF consists to determine the coefficients $\lambda{_i}$ and the polynomial $p$ such that, given a set of $N$ points $\xi$ and values $f_i$, $s$ satisfies
\begin{equation} \label{eq:rbf_con}
s(\xi) = f_i, \:\:\:i = 1,2,\ldots N
\end{equation}
If $\{p_1,\ldots, p_l\}$ is a monomial basis for polynomials of the degree of $p$, and $\c = (c_1, \ldots, c_l)\transp$ the coefficients of $p(\x)$ in this basis, then the interpolation conditions Equation~(\ref{eq:rbf_con}) can be expressed as a System of Linear Algebraic Equations (SLAE) in the form
$$
\left( \begin{array}{cc} \Phi^{\mbox{\ \ \,}} & \P \\
\P\transp & \mathbf{0} \end{array} \right)\left( \begin{array}{c}
\lamb \\ \c \end{array} \right)=\left( \begin{array}{c}  \f \\
\mathbf{0} \end{array} \right),
$$
where
$\Phi_{i,j} = \phi(\|\xi - \xj\|)$, $i = 1,\ldots,N$, $j = 1,\ldots,N$, $P_{i,j} = p_j(\xi)$, $i = 1,\ldots,N$, $j = 1,\ldots,l$
which can be simplified to
\begin{equation} \label{eq:it_sys}
\A \chi = \b
\end{equation}
where $\chi = (\lambda, \c)^T$ is the solution of the SLAE and $\b = (\f, 0)^T$ the values to be interpolated, padded with zeros.

\section{Iterative solution of CSRBF interpolation} \label{sec:proposed_method}

As mentionned previously, solving the linear system~(\ref{eq:it_sys}), is the main time consuming part of the image reconstruction process. Direct methods, similar to the one used in~\cite{Savchenko95,Turk-O'Brien99,Morse01,Wendland02}, usually fail when the size of input data exceeds a few thousands of points.

With iterative methods~\cite{Beatson99,Beatson00,Carr01,Ohtake03,Ohtake03a,Hon03}, large data sets can be addressed, although convergence is often difficult to reach, due to the conditioning of the system. An efficient way to get rid of this limitation is to increase the robustness of the algorithm by means of preconditioning techniques~\cite{Beatson99,Hon03,Saad01}.

Hybrid methods, like the non-overlapping Schwarz domain decomposition method~\cite{Unser02} and the multigrid methods~\cite{Axelsson94} have also been used. These algorithms offer powerful tools for the efficient solution of the interpolation problem, apart from the fact that their implementation in existing software requires a quite high degree of skills.

In the following a simple iterative method with a coarse space correction issued from domain decomposition methods is proposed to solve the linear system~(\ref{eq:it_sys}).
This approach
consists of a coarse space correction~\cite{2009arXiv0911.5725M,Widlund2008} applied to the solution of the interface problem arising from the domain decomposition method. In references~\cite{magoules:journal-auth:11,magoules:journal-auth:25} for graphic analysis this approach is applied directly to the solution of the linear system~(\ref{eq:it_sys}).
Each iteration of the algorithm involves a projection of the residual on a coarse space basis. With a suitable coarse space, this projection accelerates the convergence of the iterative method.
\\

The GCR (Generalized Conjugate Residual) algorithm is here considered. Not only does the GCR present similar convergence properties than the GMRES (Generalized Minimal RESidual) but it is also easier to implement in an existing software, albeit that adds some extra computation. The GCR algorithm for solving the system $\A \chi = \b$ can be written as:
\begin{algorithm}
\label{alg2}
\begin{algorithmic}[1]
 \STATE {Initialize}
\begin{displaymath}
\begin{array}{rllrll}
   \r_0 & = & \b - \A \chi_0 ; &
   \w_0 & = & \r_0 \\
\end{array}
\end{displaymath}
 \STATE {Iterate $k=0,1,2,...$ until convergence}
\begin{displaymath}
\begin{array}{rll}
   \alpha_{k} & = & \frac {(\r_k,\A \w_k)}{(\A \w_k,\A \w_k)} \\
   \chi_{k+1} & = & \chi_{k} + \alpha_k ~ \w_k \\
     \r_{k+1} & = & \r_{k} - \alpha_k  ~ \A \w_k \\
   \beta_{ik} & = & - \frac{(\A \r_{k+1},\A \w_i)}{(\A \w_i,\A \w_i)}, \mbox{ for } i= 0,1,\ldots,k \\
     \w_{k+1} & = & \r_{k+1} + \sum_{i=0}^{k} \beta_{ik}  ~ \w_i \\
\end{array}
\end{displaymath}
\end{algorithmic}
\end{algorithm}

\noindent where $k$ denotes the iteration number, $\chi_k$ the approximate solution, $\r_k = \b - \A \chi_k$ the residual vector, and $\w_k$ the search direction. 

As it is clear from this algorithm, each iteration requires a matrix-vector product, dot products and linear
combination of vectors; the matrix-vector product representing the most expensive task.
\\

The proposed method consists in projecting at each iteration, the system~(\ref{eq:it_sys}) onto a proper coarse space, this projection involving the solution of a small additional problem.
Let $\r_k$ denotes the $k$-th GCR residual, that is
\begin{equation}\label{eq:r}
 \r_k = \b - \A \chi_{k}
\end{equation}
The GCR algorithm can converge faster if, at each iteration, $\r_k$ is made orthogonal to a subspace represented by a matrix $\Q$, that is
\begin{equation}\label{eq:Qr}
  \Q^T \r_k = 0
\end{equation}
Indeed, if $\A$ is symmetric, then this weighted weak form of $\r_k = 0$ will reduce the error $\r_k$ and thus accelerate the convergence. For instance, a matrix $\Q$ with $N+4$ linearly independent columns makes the GCR method equipped with a coarse space correction converge in one iteration. Yet, it might be reminded that the subspace represented by the matrix $\Q$ should be coarse enough. Otherwise, the process of enforcing $\Q^T \r_k=0$ introduces a high unnecessary overhead.
Enforcing $\Q^T \r_k = 0$ at each GCR iteration can be achieved by means of a vector of the form $\mu = \Q \gamma$, where $\gamma$ is a vector of additional unknowns. Precisely, the vector $\chi_k$ of a GCR iteration will be replaced by a vector $\tilde \chi_k$ as follows
\begin{equation}\label{eq:chitilde}
\tilde \chi_k = \chi_k + \mu_k = \chi_k + \Q \gamma_k
\end{equation}
Then, the correction term $\mu_k = \Q \gamma_k$ enforces exactly at each iteration $k$ the optional admissible constraint $\Q^T \r_k = 0$.
Substituting Equation~(\ref{eq:chitilde}) into Equation~(\ref{eq:r}) and Equation~(\ref{eq:Qr}) shows out a projection of the initial problem  Equation~(\ref{eq:it_sys}) onto the subspace represented by $\Q$; this new problem called ``second-level'' CSRBF interpolation problem is given by:
\begin{equation}\label{eq:Qgamma}
\Q^T \A \Q \gamma_k = \Q^T (\b - \A \chi_k)
\end{equation}
From Equation~(\ref{eq:chitilde}) and Equation~(\ref{eq:Qgamma}), it follows that $\tilde \chi_k$ can be computed as
\begin{equation}\label{eq:xx0xk}
\tilde \chi_k = \chi_0 + P \chi_k
\end{equation}
where $P$ is the projector given by
$
  P = I - \Q (\Q^T \A \Q)^{-1} \Q^T \A
$
and $\chi_0$ is the initial vector given by
$
  \chi_0 = \Q (\Q^T \A \Q)^{-1} \Q^T \b
$
Finally, by substituting $\chi$ in Equation~(\ref{eq:xx0xk}) by $\tilde \chi$ in Equation~(\ref{eq:it_sys}) and multiplying the result by $P^T$, we replace the original CSRBF interpolation problem by the alternative problem
$
  P^T \A P \chi = P^T \b
$
The whole process is summarized in the following algorithm, where matrix-vector products surrounded by parentheses are simple vector variables and {\it not} actual computation. If so, only one projection of the form $P s$ and one matrix-vector product are performed per iteration.
\begin{algorithm}
\begin{algorithmic}[1]
 \STATE {Initialize}
\begin{displaymath}
\begin{array}{rllrll}
   \chi_{0} & = & \Q [\Q^T \A \Q]^{-1} \Q^T \b &&&\\
     \r_{0} & = & \b - \A \chi_0 ; & \y_0 & = & P \r_0 \\
       \w_0 & = & \y_0 , & (\A \w)_0 & = & \A \w_0 \\
\end{array}
\end{displaymath}
 \STATE {Iterate $k=1,2,...$ until convergence}
\begin{displaymath}
\begin{array}{rll}
   \zeta_{k} & = & {((\A \w)_{{k-1}}, \r_{k-1}) \over ((\A \w)_{{k-1}}, (\A \w)_{k-1})} \\
      \chi_k & = & \chi_{k-1} + \zeta_{k} ~ \w_{k-1} \\
      \r_{k} & = & \r_{k-1} - \zeta_{k} ~ (\A \w)_{k-1} \\
        \y_k & = & P \r_k \\
   (\A \y)_k & = & \A \y_k \\
        \w_k & = & \y_k - \sum_{i=0}^{i=k-1}
                   {((\A \w)_{i}, (\A \y)_k) \over ((\A \w)_{i}, (\A \w)_i)}
                    ~ \w_i \\
   (\A \w)_k & = & (\A \y)_k - \sum_{i=0}^{i=k-1}
                   {((\A \w)_{i}, (\A \y)_k) \over ((\A \w)_{i}, (\A
                   \w)_i)} ~
                   (\A \w)_i \\
\end{array}
\end{displaymath}
\end{algorithmic}
\end{algorithm}

\section{Coarse Space construction}

In~\cite{2009arXiv0911.5725M}, Mandel and Soused{\'{\i}}k  explain the principles of the design of a coarse space in a simplified way. In~\cite{Widlund2008}, Widlund shows a historically complete presentation about the development of coarse spaces for domain decomposition algorithms.
The efficiency of the coarse space correction is closely related to its key ingredient: the choice of an appropriate coarse space. Our goal here is to build such a coarse space in the context of image interpolation problem with CSRBF.

The first tentative of coarse space correction to solve CSRBF interpolation problem has been presented in~\cite{magoules:journal-auth:11}. Choosing as a coarse space basis the eigenvectors of the CSRBF interpolation problem definitely improves the convergence of the iterative method. Only few eigenvectors associated with \emph{clustered eigenvalues} are enough to accelerate the convergence of the algorithm. Unfortunately, such a choice can not be done in practice since these exact eigenvectors are too expensive to compute. Thus, a first idea is to approximate these eigenvectors numerically. An another idea presented in~\cite{magoules:journal-auth:25} consists of choosing as a coarse space basis some particular RBF. These RBF are chosen upon the RBF as the minimum set of functions required to reconstruct some basic black and white images. This choice leads to a better convergence of the iterative algorithm with the coarse space correction.
A more efficient preconditioning for the CSRBF interpolation problem was also presented in~\cite{magoules:journal-auth:25}. The authors reconstructed simple images considering as coarse space basis functions square waveforms with different frequencies and some radial basis functions with a bigger radius joined to the basis of the linear polynomial $p(\x)$. For more complex images, Daubechies wavelet basis (D4) was used.
However, as the authors noted, neither the eigenvectors associated to \emph{non-clustered eigenvalues} of the RBF interpolation problem, neither radial basis functions seems to be efficient, and the choice of a ``good'' coarse space is still an open issue.

In the following, several basis functions including trigonometric, exponentials and polynomials are investigated for CSRBF-based image reconstruction.
\\

\begin{figure}
\centering
\includegraphics[width=1.8in]{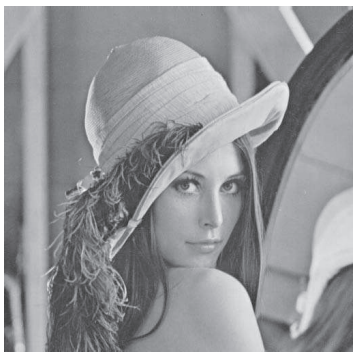}
\caption{Lena image ($512 \times 512$) used for the test case.}
\label{fig:lena}
\end{figure}

The Lena image with $262.144$ pixels, displayed Figure~\ref{fig:lena}, is used to compare the efficiency of the coarse space basis.
Former experiments applied to acoustic scattering problem~\cite{magoules:journal-auth:29} have shown the good convergence properties of the algorithm with a coarse space composed with trigonometric functions. Besides by the fourier analysis they can represent the dominant frequencies very accurately in the solution, and thus improve the convergence of the algorithm. Table~\ref{tab:coarse} shows the number of iterations requiered by the GCR with coarse space correction based on such trigonometric functions for two different stopping criteria.
The best results with the trigonometric functions are obtained with the $sinc$ and $exp$ functions.
\begin{table}[htb]
\begin{center}
\begin{tabular}{|l|l|l|l|}
\hline
\# coarse& Initial start & \# iterations ($10^{3}$) & \# iterations ($10^{6}$) \\
\hline
\multicolumn{4}{|c|}{\textsc{cosine basis}}\\
\hline
$0$ & $1$ & $88$ & $174$\\
$2$ & $0.7627$ & $84$ & $169$\\
$4$ & $0.7307$ & $83$ & $168$\\
$8$ & $0.3172$ & $82$ & $167$\\
$16$ & $0.2613$ & $78$ & $160$\\
\hline
\multicolumn{4}{|c|}{\textsc{sine basis}}\\
\hline
$0$ & $1$ & $88$ & $174$\\
$2$ & $0.7138$ & $88$ & $173$\\
$4$ & $0.7123$ & $84$ & $165$\\
$8$ & $0.7096$ & $82$ & $161$\\
$16$ & $0.7022$ & $80$ & $158$\\
\hline
\multicolumn{4}{|c|}{\textsc{tangent basis}}\\
\hline
$0$ & $1$ & $88$ & $174$\\
$4$ & $0.9832$ & $87$ & $174$\\ 
$4$ & $0.9432$ & $87$ & $173$\\
$8$ & $0.9306$ & $83$ & $170$\\
$16$ & $0.7716$ & $78$ & $152$\\
\hline
\multicolumn{4}{|c|}{\textsc{sinc basis}}\\
\hline
$0$ & $1$ & $88$ & $174$\\
$2$ & $0.7300$ & $84$ & $168$\\
$4$ & $0.7016$ & $82$ & $164$\\
$8$ & $0.2901$ & $79$ & $154$\\
$16$ & $0.2565$ & $77$ & $148$\\
\hline
\multicolumn{4}{|c|}{\textsc{exponential basis}}\\
\hline
$0$ & $1$ & $88$ & $174$\\
$2$ & $0.7679$ & $85$ & $172$\\
$4$ & $0.7581$ & $83$ & $166$\\
$8$ & $0.6348$ & $80$ & $158$\\
$16$ & $0.2625$ & $74$ & $135$\\
\hline
\multicolumn{4}{|c|}{\textsc{Gaussian basis}}\\
\hline
$0$ & $1$ & $88$ & $174$\\
$2$ & $0.8286$ & $86$ & $169$\\
$4$ & $0.8187$ & $86$ & $165$\\
$8$ & $0.6348$ & $80$ & $155$\\
$16$ & $0.4757$ & $75$ & $148$\\
\hline
\multicolumn{4}{|c|}{\textsc{Chebyshev basis}}\\
\hline
$0$ & $1$ & $88$ & $174$\\
$2$ & $0.7088$ & $84$ & $170$\\
$4$ & $0.4642$ & $82$ & $169$\\
$8$ & $0.2652$ & $79$ & $152$\\
$16$ & $0.1642$ & $69$ & $121$\\
\hline
\end{tabular}
\end{center}
\caption{Coarse space correction}
\label{tab:coarse}
\end{table}
As explained previously, the iterative method with coarse space correction converges quickly when composed of orthogonal search direction vectors. For this reason, Gaussian functions and Tchebychev functions are considered and the results reported in Table~\ref{tab:coarse}. 
Despite evaluating these functions represents almost the same computational cost than the trigonometric functions, these coarse space basis functions outperform the later one.

\section{Conclusions}
\label{sec:concluc}

In this paper, a coarse space correction is presented to solve iteratively the Radial Basis Functions interpolation problem. The method consists of an iterative method involving at each iteration a projection of the residual onto a suitable coarse space.
Numerical results illustrate the convergence properties of the proposed method for different coarse space basis for image reconstruction.

\bibliography{ref}
\bibliographystyle{abbrv}

\end{document}